\documentclass[11pt,fleqn]{article}
\usepackage{etex}
\usepackage{amssymb,latexsym,amsmath,amsfonts,graphicx, ebezier, enumitem}
\usepackage{hyperref}
\usepackage{pictex,color}
\usepackage{pict2e}
\usepackage[textsize=footnotesize]{todonotes}
    \advance\textheight by \topskip

\usepackage{tikz}
\usepackage[toc,page]{appendix}

    \numberwithin{equation}{section}

\usepackage{framed}

\makeatletter
\def\eqalign#1{\null\vcenter{\def\\{\cr}\openup\jot\m@th
  \ialign{\strut$\displaystyle{##}$\hfil&$\displaystyle{{}##}$\hfil
      \crcr#1\crcr}}\,}
\makeatother

\def\XXint#1#2#3{{\setbox0=\hbox{$#1{#2#3}{\int}$}
     \vcenter{\hbox{$#2#3$}}\kern-.5\wd0}}

\usepackage{color}

\def\beq{\begin{equation}}
\def\eeq{\end{equation}}

\newcommand{\be}{\begin{equation}}
\newcommand{\ee}{\end{equation}}

\renewcommand{\(}{\left(}
\renewcommand{\)}{\right)}

\newcommand{\e}{\mathrm{e}}

\renewcommand{\d}{\mathrm{d}}

    \def\P2n{{\rm P}_{{\rm II}}^{(n)}}

    \newtheorem{theorem}{Theorem}[section]
    
    \newtheorem{claim}{Claim}

    \newtheorem{Definition}[theorem]{Definition}
    
    \newtheorem{Remark}[theorem]{Remark}
    
    \newtheorem{Example}[theorem]{Example}
    
    \newtheorem{Assumptions}[theorem]{Assumptions}

    {\rm \trivlist \item[\hskip \labelsep{\bf Proof. }]}%
    {\hspace*{\fill}$\Box$\endtrivlist}
    {\rm \trivlist \item[\hskip \labelsep{\bf Proof}]}%
    {\hspace*{\fill}$\Box$\endtrivlist}

\makeatletter
\newcommand{\pushright}[1]{\ifmeasuring@#1\else\omit\hfill$\displaystyle#1$\fi\ignorespaces}
\newcommand{\pushleft}[1]{\ifmeasuring@#1\else\omit$\displaystyle#1$\hfill\fi\ignorespaces}
\makeatother

\begin{document}
\title{Question about integral of product of four Hermite polynomials integrated with squared weight}

\author{{\normalsize  Alexander  \textsc{Minakov}}\\[1mm]
\\
{\scriptsize Institut de Recherche en Math\'{e}matique et Physique (IRMP), Universit\'{e} catholique de Louvain
(UCL)}, \\{\scriptsize Chemin du Cyclotron 2, Louvain-la-Neuve, Belgium}
\\{\scriptsize oleksandr.minakov(at)uclouvain.be}
}
\date{\today}
\maketitle

\begin{abstract}
We found that the integral of four Hermite polynomials integrated with squared weight over the real line generates symmetric polynomials with a beautiful recursive property.
We pose a question whether that integral admit an explicit formula or not.
The question has a certain interpretation in the theory of Hermitian matrices.
\end{abstract}

\section{Claims and question}
The purpose of this note is to attract attention to the following problem.
Let $H_n(x), n\geq 0$ be Hermite polynomial orthogonal on the real line with the weight $\e^{-x^2}$ (a.k.a. ``physicists' Hermite polynomials'' to distinguish with ``probabilists' Hermite polynomials'' \cite{WikiHermite}), $$\int_{-\infty}^{+\infty}H_n(x)H_m(x)\e^{-x^2}\d x=\sqrt{\pi}\, 2^n n! \delta_{nm},$$ (where $\delta_{nm}$ is the Kronecker delta, equal to $\delta_{nm}=1$ if $n=m$ and equal to $\delta_{nm}=0$ otherwise) so that the leading coefficient is $H_n(x)=2^nx^n+\ldots, $ 
and several first polynomials are
\[H_0(x)=1,\quad H_1(x)=2x,\quad H_2(x)=4x^2-2,\quad H_3(x)=8x^3-12x,\quad H_4(x)=16x^4-48x^2+12.\]
Consider products of Hermite polynomials, and integrate them over the real line but not with the `native' weight $\e^{-x^2},$ but rather with squared weight $\e^{-2x^2};$ denote
\footnote{There must be no confusion between integral $H_n$ and polynomial $H_n(x),$ since we always write polynomial with argument.}
\begin{equation}\label{H}
\begin{split}
&H_n := \sqrt{\frac{2}{\pi}}\int_{-\infty}^{+\infty}H_n(x)\e^{-2x^2}\d x,
\\
&
H_{nm} := \sqrt{\frac{2}{\pi}}\int_{-\infty}^{+\infty}H_n(x)H_m(x)\e^{-2x^2}\d x,
\\
&
H_{nml} := \sqrt{\frac{2}{\pi}}\int_{-\infty}^{+\infty}H_n(x)H_m(x)H_l(x)\e^{-2x^2}\d x,
\\
&
H_{nmlk} := \sqrt{\frac{2}{\pi}}\int_{-\infty}^{+\infty}H_n(x)H_m(x)H_l(x)H_k(x)\e^{-2x^2}\d x.
\end{split}
\end{equation}
Our interest is in explicit formulas for $H_n, H_{nm}, H_{nml}, H_{nmlk}.$
We found them for products of no more than 3 polynomial, i.e. for $H_n, H_{nm}, H_{nml},$ and did not manage to do this for $H_{nmlk}.$
However, we claim that $H_{nmlk}$ generates symmetric polynomials
with very beautiful recursive property \eqref{Pk}.

\begin{claim}\label{claim1}
\begin{equation}
\label{Hnml}
\begin{array}{lcl}
(a)\quad & H_n &= (-1)^{\frac{n}{2}}(n-1)!!,\quad n\mbox{ is even},
\\
 & \ \quad &= 0,\quad n\mbox{ is odd},
\\\\
(b)\quad & H_{nm} &= (-1)^{\frac{n+m}{2}}(n+m-1)!!,\quad n+m\mbox{ is even},
\\
 & \ \quad &= 0,\quad n+m\mbox{ is odd},
\\\\
(c)\quad & H_{nml} &= (n+m-l-1)!!(n-m+l-1)!!(-n+m+l-1)!!,
\quad n+m+l \mbox{ is even},
\\
 & \ \quad &= 0,\quad n+m+l\mbox{ is odd},
\end{array}
\end{equation}
\end{claim}
Here $n!!$ is the double factorial, defined by $0!!=1,\ 1!!:=1,$ $n!!=n\cdot(n-2)!!$. It makes sense for all non-negative integers, and also for negative odd integers, so that
\begin{equation}\label{doublefac}(-1)!! = 1,\qquad (-3)!! = -1,\qquad (-5)!! = \frac13,\qquad (2n-1)!!(-2n-1)!!=(-1)^n.
\end{equation}
Observe that in formulas \eqref{Hnml} we use double factorial only of odd integers.
Besides, each subsequent integral $H$ embeds previous ones, so that
$H_{n0}=H_n,$ $H_{nm0}=H_{nm}.$

Dear reader, looking at formulas \eqref{Hnml} are you ready to guess
 a  formula for $H_{nmlk}?$ It is not that easy.

\begin{claim}\label{claim2}
\begin{equation}\label{HP}
H_{nmlk}
=
(-1)^k (n+m-l-k-1)!! (n-m+l-k-1)!! (-n+m+l-k-1)!!
P_k(n,m,l),
\end{equation}
where $P_k(n,m,l)$ is a symmetric polynomial of degree $2k$ with respect to variables $n,m,l.$
It is determined recursively by conditions $P_0(n,m,l)=1,$
\begin{equation}\label{Pk}
\begin{cases}
P_k(k-1,m,n) = (m-n)^2P_{k-1}(k,m,n),\\
P_k(k-2,m,n)=(m-n-1)^2(m-n+1)^2P_{k-2}(k,m,n),\\
P_k(k-3,m,n)=(m-n-2)^2(m-n)^2(m-n+2)^2P_{k-3}(k,m,n),\\
P_k(k-4,m,n)=(m-n-3)^2(m-n-1)^2(m-n+1)^2(m-n+3)^2P_{k-4}(k,m,n),\\
\ldots\\
P_k(1,m,n)=(m-n-k+2)^2(m-n-k+4)^2\cdot\ldots\cdot(m-n+k-4)^2(m-n+k-2)^2P_{1}(k,m,n),\\
P_k(0,m,n)=(m-n-k+1)^2(m-n-k+3)^2\cdot\ldots\cdot(m-n+k-3)^2(m-n+k-1)^2P_{0}(k,m,n).\\
\end{cases}
\end{equation}
or, in a shortened form, for $j=0,1,2,\ldots,k-1,$
\begin{equation}\label{Pkshort}
P_k(j,m,n)
=
\prod\limits_{l=0}^{k-j-1}(m-n-k+j+1+2l)^2\cdot P_{j}(k,m,n).
\end{equation}
\end{claim}

\subsection*{Question:} Is there a nice closed formula for $H_{nmlk}?$

\section{Discussion about $P_k(n,m,l).$}
Without going into details of proof of Claim \ref{claim2}, let us show how recursive formulas \eqref{Pk} work in practice.

\textbf{1.} For $k=0$ we have $P_0(n,m,l)=1.$

\textbf{2.} For $k=1$ formulas \eqref{Pk} read as
\[
\begin{split}
P_1(0,m,n) = (m-n)^2P_0(1,m,n) = (m-n)^2 = m^2+n^2-2mn.
\end{split}
\]
However, $P_{1}(l,m,n)$ is a symmetric w.r.t. all variables, and this allow only one option for it, namely 
\[P_1(l,m,n) = m^2+n^2+l^2-2mn-2ml-2nl.\]

\textbf{3.} For $k=2$ formulas \eqref{Pk} read as
\[
\begin{split}
&
P_2(1,m,n) = (m-n)^2P_1(2,m,n) = (m-n)^2(m^2+n^2+4-2mn-4m-4n),
\\
&
P_2(0,m,n) = (m-n-1)^2(m-n+1)^2P_0(2,m,n) = 
\(m^2+n^2-2mn-1\)^2.
\end{split}
\]
In view of symmetricity of $P_2(l,m,n)$ and the fact that degree of $P_2$ is equal to 4, the second formula allows us to determine\footnote{They are the only symmetric terms that depend on $l$ multiplicatively and have degree less or equal than 4.} the $P_2(l,m,n)$ up to terms $lmn$ and $lmn(l+m+n)$, which we can not see for $l=0:$
\[P_2(l,m,n) = \(m^2+n^2+l^2-2mn-2ml-2nl-1\)^2+C_1mnl+C_2 mnl(m+n+l),\]
and substituting the latter into the formula for $P_2(1,m,n)$ allows us to determine the constants $C_1,$ $C_2:$ we find $C_1=-16$, $C_2=0.$

\textbf{4.} We can continue for bigger values of $k$ in the same manner as for $k=1,2.$ This gives us {subsequently}~\footnote{We repeat signs after line breaking, so `-' at the end of a line combined with `-' in the beginning of the next line means~`-'.}
\[
\begin{split}
P_2(m,n,l) =& [m^2+n^2+l^2-2(mn+ml+nl)-1]^2-16mnl = \ (P_1-1)^2-16mnl,
\\
P_3(m,n,l) =& \left(P_1^3-8P_1^2+(-48m n l+16)P_1-64m n l(m+n+l-3)\right),
\\
P_4(m,n,l) =& \left(P_1^4-20P_1^3+(-96m n l+118)P_1^2+(960m n l-180)P_1-\right.
\\
&
-256m n l(m^3+n^3+l^3)+256mnl(m^2l+m^2n+l^2m+l^2n+n^2l+n^2m)
+
\\
&
\left.
+2304m^2n^2l^2-1536mnl(mn+ml+ln)+2304mnl(m+n+l)-3936mnl+81\right),
\\\\
P_5(m,n,l) =& P_1^5\ -\ 160mnlP_1^3\ -\ 40P_1^4\ -\ 640mnl(m+n+l)P_1^2\ +\ 528P_1^3\ +
\\
&
+3840m^2n^2l^2P_1 \ + \ 3200mnlP_1^2 - 7680mnl(m^3(n+l)+n^3(m+l)+l^3(m+n))+
\\
&
+14336mnl(m^3+n^3+l^3)
\ +\ 15360mnl(m^2n^2+n^2l^2+m^2l^2)\ -\ 2560P_1^2+
\\
&
+ 33280m^2n^2l^2(m+n+l)\  \ -
\ 26624mnl(m^2(l+n)+l^2(m+n)+n^2(m+l))-
\\
&
- 17920mnlP_1\ -\ 172032m^2n^2l^2
\ +\ 122880mnl(mn+nl+ml) \ - 
\\
&
- 102400mnl(m+n+l)\ +\ 4096 P_1\ +\ 122880 mnl.
\end{split}
\]
Here we denoted for shortness $P_1=P_1(l,m,n).$

\textbf{4.} It is possible to continue further, however, the expressions for the first $P_k$ do not suggest any nice structure for the polynomials of higher degree.

\section{Discussion about proofs}\label{sect_proof}
Since we did not reach our goal, which was an explicit expression for $H_{nmlk},$ I do not think that complete proofs of Claims \ref{claim1}, \ref{claim2} will engage a reader.
However, I would like to offer a discussion about possible ways to prove them and about difficulties on those ways.

Hermite polynomials admit a generating function \cite{WikiHermite},
\[e^{2xr-r^2}=\sum\limits_{n=0}^{\infty}\frac{H_n(x)}{n!}r^n.\]
Taking product of $\e^{2xr-r^2},$ $\e^{2xs-s^2},$ and so on, multiplying by $\e^{-2x^2}$ and integrating, we find subsequently

\begin{equation}\label{prod_H_1}\sum\limits_{n=0}^{\infty}\frac{r^n}{n!}\int\limits_{-\infty}^{+\infty}H_n(x) \ e^{-2x^2}dx =\int\limits_{-\infty}^{+\infty}e^{2xr-r^2}e^{-2x^2}dx=e^{\frac{-1}{2}r^2}\sqrt{\frac{\pi}{2}}.
\end{equation}

\begin{equation}\label{prod_H_2}\sum\limits_{m,n=0}^{\infty}\frac{r^ns^m}{n!m!}\int\limits_{-\infty}^{+\infty}H_n(x)H_m(x) \ e^{-2x^2}dx =\int\limits_{-\infty}^{+\infty}e^{2xr-r^2+2xs-s^2}e^{-2x^2}dx=e^{\frac{-1}{2}(r-s)^2}\sqrt{\frac{\pi}{2}}.
\end{equation}
\begin{equation}\label{prod_H_3}
\begin{split}
&\sum\limits_{m,n,l=0}^{\infty}\frac{r^ns^mt^l}{n!m!l!}\int\limits_{-\infty}^{+\infty}H_n(x)H_m(x)H_l(x) \ e^{-2x^2}dx =
\\
&\hskip4cm 
=\int\limits_{-\infty}^{+\infty}e^{2xr-r^2+2xs-s^2+2xt-t^2}e^{-2x^2}dx=e^{-\frac12(r-s-t)^2+2st}\sqrt{\frac{\pi}{2}},
\end{split}
\end{equation}
\begin{equation}\label{prod_H_4}
\begin{split}
&
\sum\limits_{m,n,l,j=0}^{\infty}\frac{r^ns^mt^lu^j}{n!m!l!j!} \int\limits_{-\infty}^{+\infty}H_n(x)H_m(x)H_l(x)H_j(x) 
\e^{-2x^2}dx 
=
\\
&
\hskip1cm =\int\limits_{-\infty}^{+\infty}e^{2xr-r^2+2xs-s^2+2xt-t^2+2xu-u^2}e^{-2x^2}dx=e^{\frac{-1}{2}(r-s-t-u)^2+2(st+su+ut)}\sqrt{\frac{\pi}{2}}.
\end{split}
\end{equation}
Exponent terms in the r.h.s. are written not in a symmetric way, but rather in a way most convenient for expanding them in power series of the variables $r,s,t,u,$ where applicable.

Expanding, for \eqref{prod_H_1} and \eqref{prod_H_2} we obtain
\begin{equation}\label{prod1expanded}
\sum\limits_{n=0}^{\infty}\frac{H_{n} r^n}{n!}
=\sum\limits_{\tiny{
\begin{array}{ccc}n\geq0,
\\
n \mbox{ even}
\end{array}}}
\frac{(-1)^{\frac{n}{2}}C_{n}^n r^n}{(\frac{n}{2})!\,2^{\frac{n}{2}}},
\end{equation}
\begin{equation}\label{prod2expanded}
\sum\limits_{m,n=0}^{\infty}\frac{H_{nm} r^ns^m}{n!m!} 
=\sum\limits_{\tiny{
\begin{array}{ccc}m\geq0, n\geq0,
\\
m+n \mbox{ even}
\end{array}}}
\frac{(-1)^{\frac{m-n}{2}}C_{n+m}^n r^n s^m}{(\frac{n+m}{2})!\,2^{\frac{n+m}{2}}},
\end{equation}
where $C_p^k=\frac{p!}{k!(p-k)!}$, $k=0,1,\ldots,p$ is a binomial coefficient. After elementary transformations these formulas give us statements (a), (b) of Claim \ref{claim1}.
It is important that in both formulas \eqref{prod1expanded},
\eqref{prod2expanded} the number of indices of summation in the r.h.s. coincides with the number of indices in the l.h.s.

However, this pattern breaks and things start to getting worse already for \eqref{prod_H_3}. 
Namely, \eqref{prod_H_3} gives rise to a power series with 4 indices
\footnote{Indeed, exponent gives 1 index, and sum $(r-s-t)^{2k}$ gives 3 indices.}, which is more than 3 indices on the l.h.s.
Situation is even worse for formula \eqref{prod_H_4}, which produces 
8 indices of summation \footnote{Indeed, 1 index is produced by exponent, 4 indices are produced by $(r-s-t-u)^{2k}$ term, and 3 indices are produced by $(st+su+ut)^{n-2k}$ term.}.
Hence, it is not feasible to simplify summation formulas directly,
but since the expressions for these formulas are already found in Claims \ref{claim1}, \ref{claim2}, one can prove them by induction.

\section{Discussion about $P_k(n,m,l),$ continued.}
It is only formula \eqref{HP} which is not proved in Claim 2.
Once we assume that \eqref{HP} holds true, the recurrence formulas \eqref{Pk} for polynomials $P_k(n,m,l)$ follow from symmetricity of $H_{k,l,m,n}$ and the triviality formula
\begin{equation}\label{trivial}\begin{split}
&(l+m-n-k-1)!!(l-m+n-k-1)!!(-l+m+n-k-1)!!\cdot
\\
&\cdot(-l+m-n+k-1)!!(-l-m+n+k-1)!!(l-m-n+k-1)!!=1.
\end{split}
\end{equation}

\section{Some further formulas}
Formula \eqref{Pkshort} can be rewritten in the form
\[\frac{P_{k}(l,m,n)}{P_{l}(k,m,n)}=\frac{(m-n+k-l-1)!!^2}{(m-n-k+l-1)!!^2}.\]

Furthermore, from \eqref{HP} and \eqref{trivial} it follows that
\[\begin{split}
H^2_{klmn}&=(-1)^{\frac{l+m+n+k}{2}}\((-l+m+n-k-1)!!\)^2P_l(k,m,n)P_k(l,m,n)=
\\
 &=
\frac{(-1)^{\frac{l+m+n+k}{2}}}{\((l-m-n+k-1)!!\)^2}P_l(k,m,n)P_k(l,m,n).
\end{split}
\]
In the above formula, changing $k, l $ with $m , n$ and then multiplying, we get (using the last formula in \eqref{doublefac})
\[H^4_{klmn}=P_{l}(k,m,n)P_{m}(k,n,l)P_{n}(k,m,l)P_{k}(l,m,n)(-1)^{l+m+n+k}.\]

\section{Recurrent formula for $H_{nmlk}$.}
The next idea belongs to \textbf{ Di Yang }\cite{DiYang}:
using the recurrence formula for Hermite polynomials
\[H_n(x)=2xH_{n-1}(x)-2(n-1) H_{n-2}(x),\qquad H'_n(x)=2n H_{n-1}(x),\]
we obtain
\[
\begin{split}
&H_{nmlj} = \sqrt{\frac{2}{\pi}}\int\limits_{-\infty}^{+\infty}\(2xH_{n-1}(x)-2(n-1)H_{n-2}(x)\)H_m(x)H_l(x)H_j(x)\e^{-2x^2}\d x=
\\
&=-\frac12\sqrt{\frac{2}{\pi}}\int\limits_{-\infty}^{+\infty}\(H_{n-1}(x)\)H_m(x)H_l(x)H_j(x)\(\e^{-2x^2}\)'\d x-2(n-1)H_{n-2,m,l,j}=
\\
&=
\frac12 \left[2(n-1)H_{n-2,m,l,j}+2mH_{n-1,m-1,l,j}+2lH_{n-1,m,l-1,j}+
	2jH_{n-1,m,l,j-1}\right]-2(n-1)H_{n-2,m,l,j}=
\\
&=-(n-1)H_{n-2,m,l,j}+mH_{n-1,m-1,l,j}+lH_{n-1,m,l-1,j}
	+jH_{n-1,m,l,j-1}.
\end{split}
\]
This gives us another way to obtain $H_{nmlj}$
\[
H_{nmlj}=-(n-1)H_{n-2,m,l,j}+mH_{n-1,m-1,l,j}+
lH_{n-1,m,l-1,j}+jH_{n-1,m,l,j-1}.
\]
For polynomials $P_n(m,l,j)$ this relation becomes
\[\begin{split}
P_{n}(m,l,j)=&
-(n-1)(m+l-j-n+1)(m-l+j-n+1)(-m+l+j-n+1)P_{n-2}(m,l,j)-
\\
&-m(-m+l+j-n+1)P_{n-1}(m-1,l,j)-
\\
&
-l(m-l+j-n+1)P_{n-1}(m,l-1,j)-
\\
&
-j(m+l-j-n+1)P_{n-1}(m,l,j-1).
\end{split}
\]

This formula might be used to prove Claim \ref{claim2}, but it still does not allow to obtain or guess an explicit expression for $H_{nmlj}.$

\section{The original motivation.}
The study originated from the interest of \textbf{Khazhgali Kozhasov}\cite{Khazh} in the following integral of determinant:
\[
D_n = \int\limits_{-\infty}^{+\infty}
\det\left[H_{n+i+j}(x)\right]_{i,j=0}^{3}\,\e^{-2x^2} \d x
\]
or at least in asymptotics of $D_n$ for large $n.$ The $D_n$ is related to the volume of Hermitian matrices with two coinciding eigenvalues \cite{Khazh}.

Expanding the above determinant, we will obtain a linear combination of terms $H_{nmlk},$ where all $n,m,l,k$ differ from each other by no more than 6.

It is \textbf{Giulio Ruzza's}\cite{Giulio} idea to take a look at large $n$ asymptotics of $H_n$, $H_{nm}$, $H_{nml}$ (with $n,m,l$ being close to each other) and try to guess the asymptotic formula for $H_{nmlk}$.
On this way one gets, as $n\to\infty,$
\[
H_n\sim\e^{\frac{n}{2}\ln\frac{n}{\e}}\sqrt{2}(-1)^{\frac{n}{2}}\(1+\mathcal{O}(\frac{1}{n})\),
\]
\[
H_{nn}\sim\e^{n\ln\frac{n}{\e}+n\ln2}\sqrt{2}\(1+\mathcal{O}(\frac{1}{n})\),
\]
\[
H_{nnn}\sim\e^{\frac{3n}{2}\ln\frac{n}{\e}}2\sqrt{2}\(1+\mathcal{O}(\frac{1}{n})\),
\]
\[
\frac{H_{nnn}}{H_{n+1,n-1,n}}\sim1-\frac{2}{n}+\mathcal{O}(n^{-2}),\quad \frac{H_{n,n,n}}{H_{n-1,n-1,n}}\sim n-1,\quad 
\frac{H_{n,n,n}}{H_{n-2,n,n}}\sim n-3+\mathcal{O}(n^{-1}),
\]
\[
\frac{H_{n,n,0}}{H_{n+1,n-1,0}}\sim-1,\quad \frac{H_{n,n,0}}{H_{n-1,n-1,0}}\sim 2n-1,
\]
and it is possible to write down asymptotics of other quantities, in order to get an intuition, but we do not pursue this idea further.

\end{document}